\documentclass[11pt, a4paper]{article}

\usepackage{array}

\newcolumntype{"}{@{\hskip\tabcolsep\vrule width 1pt\hskip\tabcolsep}}

\makeatother

\usepackage{epsfig}
\usepackage{subfigure}
\usepackage{amssymb}
\usepackage{amsthm}
\usepackage{mathrsfs}
\usepackage{hhline}
\usepackage{longtable}
\usepackage{amsmath}
\usepackage{float}
\usepackage{multicol}
\usepackage{multirow}
\usepackage{cite}
\usepackage{enumerate}
\usepackage{xcolor}
\usepackage{caption}
\usepackage{tabu}

\def\ms{\medskip}
\def\nt{\noindent}

\def\purple{\color{purple}}
\definecolor{vividviolet}{rgb}{0.62, 0.0, 1.0}
\def\rsq{\hspace*{\fill}$\Box$\medskip}

\def\Z{\mathbb Z}
\newtheoremstyle{de}
  {10pt}          
  {10pt}  
  {\rm}  
  {}
  {\bf}  
  {. }    
  { }    
  {}     
\theoremstyle{de}

\newtheorem{example}{Example}[section]

\newtheorem{problem}{Problem}[section]

\newtheoremstyle{theorem}
  {10pt}          
  {10pt}  
  {\it}  
  {}
  {\bf}  
  {. }    
  { }    
  {}     
\theoremstyle{theorem}
\newtheorem{theorem}{Theorem}[section]

\newtheorem{conjecture}{Conjecture}[section]

\numberwithin{equation}{section}
\def\Z{\mathbb{Z}}

\graphicspath{%
    {converted_graphics/}
    {/}
}

\begingroup\makeatletter\ifx\SetFigFont\undefined%
\gdef\SetFigFont#1#2#3#4#5{%
  \reset@font\fontsize{#1}{#2pt}%
  \fontfamily{#3}\fontseries{#4}\fontshape{#5}%
  \selectfont}%
\fi\endgroup%

\def\Lau{\color{vividviolet}}

\begin{document}
\baselineskip18truept
\normalsize
\begin{center}
{\mathversion{bold}\Large \bf Constructions of local antimagic 3-colorable graphs\\ of fixed odd size | matrix approach}

\bigskip
{\large  G.C. Lau$^{a,}$\footnote{Corresponding author.},  W.C. Shiu{$^b$}, K. Premalatha{$^c$}, M. Nalliah{$^d$}}\\

\medskip

\emph{{$^a$}77D, Jalan Subuh, }\\
\emph{85000, Segamat, Johor, Malaysia.}\\
\emph{geeclau@yahoo.com}\\

\medskip
\emph{{$^b$}Department of Mathematics,}\\
\emph{The Chinese University of Hong Kong,}\\
\emph{Shatin, Hong Kong, P.R. China.}\\
\emph{wcshiu@associate.hkbu.edu.hk}\\
\medskip

\emph{{$^c$}Department of Mathematics,}\\
\emph{Sri Shakthi Institute of Engineering and Technology,}\\
\emph{ Coimbatore, 641062 Tamil Nadu, India.}\\
\emph{premalatha.sep26@gmail.com}\\

\medskip
\emph{{$^d$}Department of Mathematics, School of Advanced Sciences, }\\
\emph{Vellore Institute of Technology, Vellore, 632014, Tamil Nadu, India.}\\
\emph{nalliah.moviri@vit.ac.in}\\
\end{center}

\begin{center} Dedicated to Prof Sin-Min Lee on the occasion of his 80th birthday \end{center}

\begin{abstract}
An edge labeling of a connected graph $G = (V, E)$ is said to be local antimagic if there is a bijection $f:E \to\{1,\ldots ,|E|\}$ such that for any pair of adjacent vertices $x$ and $y$, $f^+(x)\not= f^+(y)$, where the induced vertex label $f^+(x)= \sum f(e)$, with $e$ ranging over all the edges incident to $x$.  The local antimagic chromatic number of $G$, denoted by $\chi_{la}(G)$, is the minimum number of distinct induced vertex labels over all local antimagic labelings of $G$. In this paper, we give three ways to construct a $(3m+2)\times (2k+1)$ matrix that meets certain properties for $m=1,3$ and $k\ge 1$. Consequently, we obtained many (disconnected) graphs of size $(3m+2)(2k+1)$ with local antimagic chromatic number 3.
\ms


\end{abstract}

\baselineskip18truept
\normalsize

\section{Introduction}
Let $G=(V, E)$ be a connected graph of order $p$ and size $q$.
A bijection $f: E\rightarrow \{1, 2, \dots, q\}$ is called a \textit{local antimagic labeling}
if $f^{+}(u)\neq f^{+}(v)$ whenever $uv\in E$,
where $f^{+}(u)=\sum_{e\in E(u)}f(e)$ and $E(u)$ is the set of edges incident to $u$.
The mapping $f^{+}$ which is also denoted by $f^+_G$ is called a \textit{vertex labeling of $G$ induced by $f$}, and the labels assigned to vertices are called \textit{induced colors} under $f$.
The \textit{color number} of a local antimagic labeling $f$ is the number of distinct induced colors under $f$, denoted by $c(f)$.  Moreover, $f$ is called a \textit{local antimagic $c(f)$-coloring} and $G$ is {\it local antimagic $c(f)$-colorable}. The \textit{local antimagic chromatic number} $\chi_{la}(G)$ is defined to be the minimum number of colors taken over all colorings of $G$ induced by local antimatic labelings of $G$~\cite{Arumugam}. Let $G+H$ and $mG$ denote the disjoint union of graphs $G$ and $H$, and $m$ copies of $G$, respectively. For integers $a < b$, let $[a,b] = \{n\in\Z\;|\; a\le n\le b\}$.

\ms\nt Haslegrave~\cite{Haslegrave} proved that every connected graph except $K_2$ admits a local antimagic coloring. Thus, so does every disconnected graph with $K_2$.  In~\cite{LNS, LPSN-Proy, LS, LSN}, the authors obtained the local antimagic chromatic number for various families of connected graphs using different matrix representations techniques including matrices with repeated entries or with empty cells.  Motivated by this, we give three ways to construct  a $(3m+2)\times (2k+1)$ matrix that meets various properties. We first generalized $\chi_{la}(P_3 \vee O_m)=3$ in~\cite[Corollary 3.18]{LSN} to $\chi_{la}(nP_3 \vee O_m)=3$, for odd $n, m\ge 3$, where $O_m$ is the null graph of order $m$. Consequently, we obtained many new (disconnected) graphs with local antimagic chromatic number 3.  The results for odd number of rows and even number of column matrix that give many new results including $\chi_{la}(nP_3 \vee O_m)=3$, for even $n, m\ge 2$, will appear in another paper.

\section{$m=1$} 

We shall always refer to the following table (with entries in $[1,10k+5]$ bijectively) to get our results in this section.
\[\fontsize{8}{11}\selectfont
\begin{tabu}{|c|[1pt]c|c|c|c|c|c|[1pt]c|c|c|c|c|c|}\hline
i & 1 & 2 & 3 & \cdots & k & k+1 & k+2 & k+3 & \cdots & 2k-1 & 2k & 2k+1 \\\tabucline[1pt]{-}
f(u_iw_i) & 1 & 3 & 5 & \cdots & 2k-1 & 2k+1 & 2 & 4 & \cdots & 2k-4 & 2k-2 & 2k \\\hline
f(v_iw_i) & 3k+2 & 3k+1 & 3k & \cdots & 2k+3 & 2k+2 & 4k+2 & 4k+1 & \cdots & 3k+5 & 3k+4 & 3k+3 \\\hline
f(x_iw_i) & 6k+3 & 6k+2 & 6k+1 & \cdots & 5k+4 & 5k+3 & 5k+2 & 5k+1 & \cdots & 4k+5 & 4k+4 & 4k+3 \\\hline
f(x_iu_i) & 10k+5 & 10k+3 & 10k+1 & \cdots & 8k+7 & 8k+5 & 10k+4 & 10k+2 & \cdots & 8k+10 & 8k+8 & 8k+6 \\\hline
f(x_iv_i) & 7k+4 & 7k+5 & 7k+6 & \cdots & 8k+3 & 8k+4 & 6k+4 & 6k+5 & \cdots & 7k+1 & 7k+2 & 7k+3 \\\hline
\end{tabu}
\]

%
Observe that
\begin{enumerate}[(1)]
  \item In each column, the sum of the first 3 row entries is a constant $S_1=9k+6.$
  \item In each column, the sum of the rows 1 and 4 (respectively, 2 and 5) entries is $S_2=10k+6$.
  \item The sums of the last 3 row entries for columns 1 to $2k+1$ respectively form an arithmetic progression with first sum is $23k+12$ to last sum is $19k+12$ in decreasing of 2.
  \item Since the entries of row 3 form an arithmetic progression in decreasing of 1, by (3) the sum of entries of rows 4 and 5 form  an arithmetic progression in decreasing of 1.
  \item The sum of all the last 3 row entries is $S = (7k+4)(6k+3)$.

  \item If $2k+1=rs \ge 3$, then we can divide the table into $r\ge 3$ blocks of $s$ column(s) with the $j$-th block containing $(j-1)s+1, (j-1)s+2, \ldots,js$ columns. By (3) and (4), the sum of the row 3 entries in the $j$-th block, and the rows 4 and 5 entries in the $(r+1-j)$-th block is $S_3=s(21k+12)$. Similarly, the sum of the last three row entries in $(r+1)/2$ block is also $S_3$.
\end{enumerate}

\nt Denote $nP_3\vee O_1$ by $FB(n)$, the {\it  fan graph with $n$ blades}. Note that $FB(1)=P_3\vee O_1 \cong K_{1,2}\vee O_1$ is the fan graph $F_3$ of order 4.

\begin{theorem}\label{thm-FBn} For odd $n\ge 1$, $\chi_{la}(FB(n)) = 3$. \end{theorem}

\begin{proof} Suppose $n=2k+1\ge 1$. Thus, $FB(2k+1)$ is of size $10k+5$. In~\cite[Theorem 3.7]{LSN}, the authors proved that $\chi_{la}(FB(1)) = 3$. We now consider $k\ge 1$. Let $G=(2k+1)FB(1)$ with $V(G) = \{u_i,v_i,w_i,x_i\;|\;1\le i\le 2k+1\}$ and $E(G) = \{u_iw_i, v_iw_i, x_iw_i, x_iu_i, x_iv_i\;|\;1\le i\le 2k+1\}$. Observe that table above gives a bijective edge labeling $f$ of  $G$ using integers in $[1,10k+5]$. Merging the vertices $x_1$ to $x_{2k+1}$ to get a vertex $x$, we get $FB(n)$ with induced vertex labels $f^+(u_i) = f^+(v_i) = 10k+6$, $f^+(w_i) = 9k+6$, for $1\le i\le 2k$, and $f^+(x) = \sum\limits^{2k+1}_{i=1} f^+(x_i) =  (7k+4)(6k+3)$. Thus, $\chi_{la}(FB(n))\le 3$.  Since $\chi_{la}(FB(n))\ge \chi(FB(n))=3$, the theorem holds.
\end{proof}

\begin{theorem}\label{thm-tFBs} For odd $t,s\ge 3$, $\chi_{la}(tFB(s)) = 3$. \end{theorem}

\begin{proof} Let $ts = 2k+1\ge 9$. Consider $tsFB(1)$ with the bijective edge labeling as defined in the proof of Theorem~\ref{thm-FBn}.
By Observation (3) above, we can arrange integers in $\{f^+(x_i)\;|\; 1\le i\le 2k+1\}=\{23k+12, 23k+10, \ldots, 19k+12\}$ into an $t\times s$ magic rectangle $M$ with constant row sum $(2k+1)(21k+12)/t = s(21k+12)$. For each row $1\le a\le t$ in $M$, we merge the $x_i$'s with incident edge labels sum in row $a$ to get a new vertex $y_a$ of degree $3s$ with incident edge label sum $s(21k+12)$. We now have $t$ copies of $FB(s)$ that admits a local antimagic 3-coloring with induced vertex labels $10k+6, 9k+6, s(21k+12)$. Thus, $\chi_{la}(t(FB(s))) \le 3$. Since $\chi_{la}(t(FB(s)))\ge \chi(t(FB(s)))=3$, the theorem holds.
\end{proof}


\begin{example}\label{eg-FB9} Following is a labeling for the graph $9FB(1)$ according to the proof of Theorem~\ref{thm-tFBs}. 

\begin{figure}[H]
\centerline{\epsfig{file=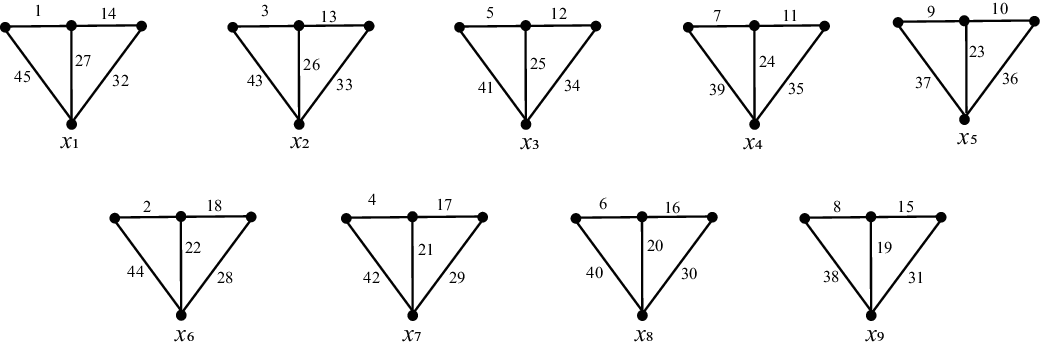, width=11cm}}
\caption{The graph $9FB(1)$ for the graphs $FB(9)$ and $3FB(3)$.}\label{fig:fb9}

\end{figure}

\nt By merging the vertices $x_1, \ldots, x_9$ of $9FB(1)$ into the single vertex $x$, we obtain a labeling for $FB(9)$.  Consider $t=s=3$, $3FB(3)$ is obtained from the graph $9FB(1)$ by merging the vertices in $\{x_1,x_5,x_9\}$, $\{x_3,x_4,x_8\}$, $\{x_2,x_6,x_7\}$. \rsq
\end{example}

\nt  For $s\ge 1$, let $S_1$ and $S_2$ be two copies of $sP_3$. Let $DF(2s)$ be the {\it diamond fan} graph obtained from $S_1+S_2$ by joining a vertex $y$ (resp., $z$) to every degree 1 vertex of $S_1$ (resp., $S_2$) and every degree 2 vertex of $S_2$ (resp., $S_1$).
For $r\ge 1$, let $DF(r, 2s) = rDF(2s) + FB(s)$. Note that $DF(r, 2s)$ is of size $5(2r+1)s$ has
\begin{enumerate}[1.]
\item $(4r+2)s$ vertices of degree 2; these vertices will be denoted by $u_i$ and $v_i$, $1\le i\le (2r+1)s$;
\item $(2r + 1)s$ vertices of degree 3; these vertices will be denoted by $w_i$, $1\le i\le (2r+1)s$;
\item $2r+1$ vertices of degree $3s$; these vertices will be denoted by $x$,  $y_i$ and $z_i$, $1\le i\le r$. Hence $s$ can be 1.
\end{enumerate}

\begin{theorem}\label{thm-nDF2sodd} For $r\ge 1$ and odd $s\ge 1$,  $\chi_{la}(DF(r, 2s))=3$.  \end{theorem}

\begin{proof} Let $ (2r+1)s = 2k+1\ge 3$. Begin with $(2r+1)s = 2k+1\ge 3$ copies of $FB(1)$ with edge labeling as in the proof of Theorem~\ref{thm-FBn}. Partition the $(2r+1)sFB(1)$ into $(2r+1)$ blocks of $s$ copies of $F(1)$ such that the $j$-th block has vertices $x_{(j-1)s+a}$ for $1\le j\le 2r+1, 1\le a\le s$.  For $j\in [1,2r+1]\setminus \{r+1\}$, split each $x_i$ into $x^1_i$ and $x^2_i$ such that $x^1_i$ is adjacent to $w_i$ and $x^2_i$ is adjacent to $u_i,v_i$. For each $1\le j\le r$ and $1\le a\le s$, merge the vertices in $\{x^1_{(j-1)s+a}, x^2_{(2r+1-j)s+a}\}$ to get vertex $y_j$, and the vertices in $\{x^2_{(j-1)s+a}, x^1_{(2r+1-j)s+a}\}$ to get vertex $z_j$. We now have $r$ copies of $DF(2s)$. The $(r+1)$-st block is used to construct the $FB(s)$ by merging vertices $x_{rs+1}$ to $x_{(r+1)s}$.

\ms\nt By Observation (6) above, it is easy to verify that the graph obtained is $DF_r(2s)$ that admits a local antimagic 3-coloring with degree 2 (respectively, $3$ and $3s$) vertices having induced vertex labels $10k+6$ (respectively, $9k+6$ and $(21k+12)s$).  Thus, $\chi_{la}(DF(r, 2s)) \le 3$.

\ms\nt Since the $FB(s)$ component is a tripartite graph, we have $\chi_{la}(DF(r, 2s))\ge \chi(DF(r, 2s))=3$. Hence We have $\chi_{la}(DF(r, 2s))=3$.
\end{proof}

\begin{example} Using $FB(9)$ by taking $r=1$ and $s=3$, we have $DF(1, 6)=DF(6)+FB(3)$ as shown below. We can get $FB(3)$ from the graphs in Example~\ref{eg-FB9} by merging vertices $x_4$, $x_5$, $x_6$.

\begin{figure}[H]
\begin{center}
\epsfig{file=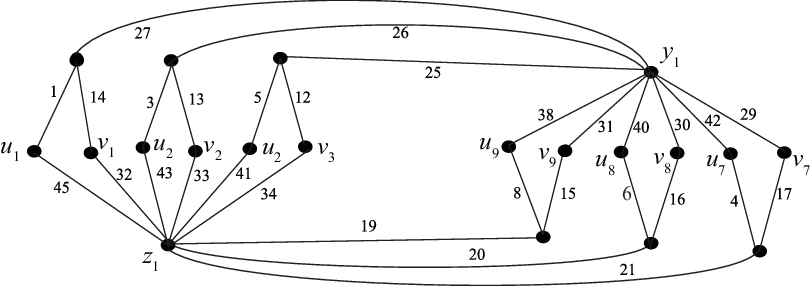, width=8cm}\qquad\epsfig{file=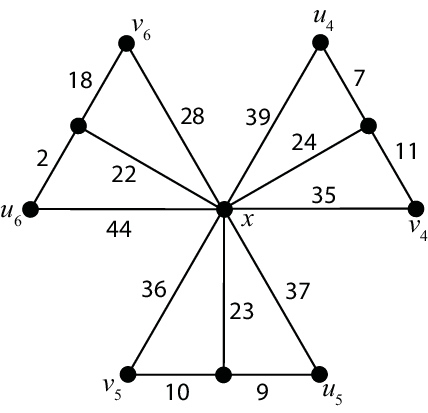, width=3cm}
\caption{Graph $DF(1, 6)=DF(6)+FB(3)$.}\label{fig:DF6}
\end{center}
\end{figure}

\nt Using $FB(9)$, we can also get $DF(4, 2) = 4DF(2) + FB(1)$. The $i$-th component of the $4DF(2)$ can be obtained from the $i$-th and the $(10-i)$-th $FB(1)$, while the $FB(1)$ component is the 5-th $FB(1)$ in Example~\ref{eg-FB9}.
\rsq
\end{example}

\nt For $r,s\ge 2$, $1\le i\le r$, let $G_i$ be the $i$-th component of $rFB(s)$ with vertex set\break $\{u_{i,j}, v_{i,j}, w_{i,j}, x_i\;|\;1\le j\le s\}$ and edge set $\{u_{i,j}w_{i,j}, v_{i,j}w_{i,j}, x_iw_{i,j}\;|\;1\le j\le s\}$. Note that $u_{i,j}$ and $v_{i,j}$ are of degree 2, $w_{i,j}$ is of degree 3, and $x_i$ is of degree $3s$.


\ms\nt Denote by $FB^1(r,s)$ the graph obtained from $G_1, \dots, G_r$ by merging the vertices in $\{u_{i,j}\;|\; 1\le i\le r\}$ and in $\{v_{i,j}\;|\; 1\le i\le r\}$ respectively for each $1\le j\le s$. Let the new vertices be $u_j$ and $v_j$ respectively for $1\le j\le s$. Note that $FB^1(r,s)$ has $rs$ vertices of degree 3, $2s$ vertices of degree $2r$, and $r$ vertices of degree $3s$.

\ms\nt Denote by $FB^2(r,s)$ the graph obtained from $G_1, \dots, G_r$ by merging the vertices in $\{w_{i,j}\;|\; 1\le i\le r\}$ respectively for $1\le j\le s$. Let the new vertices be $w_j$ for $1\le j\le s$. Note that $FB^2(r,s)$ has $2rs$ vertices of degree 2, $s$ vertices of degree $3r$, and $r$ vertices of degree $3s$.


\begin{theorem}\label{thm-FBrs12} For odd $r,s\ge 3$, $\chi_{la}(FB^1(r,s)) = 3$ if $k\not\equiv 2\pmod{4}$, and  $\chi_{la}(FB^2(r,s))=3$. \end{theorem}

\begin{proof} Let $rs = 2k+1$.  Using the $r$ copies of $FB(s)$ in Theorem~\ref{thm-tFBs} and the corresponding local antimagic 3-coloring by merging the degree 2 vertices as defined for $FB^1(r,s)$, we can immediately conclude that $FB^1(r,s)$ admits a local antimagic 3-coloring having $rs$ (respectively, $2s$ and $r$) vertices of degree 3 (respectively, $2r$ and $3s$) with induced vertex label $9k+6$ (respectively, $r(10k+6)$ and $s(21k+12)$). Suppose $r(10k+6)=s(21k+12)$. Since $s$ is odd, we get $(2k+1)(10k+6) = s^2(21k+12)$. Thus, $(2k+1)(2k+2)\equiv k \pmod{4}$ so that $k\equiv 2\pmod{4}$. Therefore, all the induced vertex labels are distinct if $k\not\equiv 2\pmod{4}$.

\ms\nt Similarly, by merging the degree 3 vertices as defined for $FB^2(r,s)$, we can immediately conclude that $FB^2(r,s)$ admits a local antimagic 3-coloring having $2rs$ (respectively, $s$ and $r$) vertices of degree 2 (respectively, $3r$ and $3s$) with induced vertex label $10k+6$ (respectively, $r(9k+6)$ and $s(21k+12)$). If  $r(9k+6)=s(21k+12)$, then $(2k+1)(9k+6)\equiv s^2(21k+12)$. Thus, $(2k+1)(k+2)\equiv k \pmod{4}$ so that $k^2\equiv -1\pmod{4}$, a contradiction.

\ms\nt Thus, $\chi_{la}(FB^1(r,s)) \le 3$ if $k\not\equiv 2\pmod{4}$, and $\chi_{la}(FB^2(r,s))\le 3$. Since $\chi_{la}(FB^i(r,s))\ge \chi(FB^i(r,s))=3$, $i=1,2$, the theorem holds.
\end{proof}

\nt Consider $DF(r, 2s)$ as in Theorem~\ref{thm-nDF2sodd} with $(2r+1)s = 2k+1\ge 3$. Let $V_1$ be the set of $(4r+2)s\ge 6$ vertices of degree 2 with induced vertex label $10k+6$. Partition $V_1$ into $2s\ge 4$ blocks of size $2r+1\ge 3$ such that each block contains exactly one degree 2 vertex of $FB(s)$ and exactly two degree 2 vertices of each $DF(2s)$ component that are without common neighbors.  Let $V_2$ be the set of $(2r+1)s\ge 9$ vertices of degree 3 with induced vertex label $9k+6$. Partition $V_2$ into $s\ge 3$ blocks of size $2r+1\ge 3$ such that each block contains exactly one degree 3 vertex of $FB(s)$ and exactly two degree 2 vertices of each $DF(2s)$ component that are without common neighbors. For $1\le i\le 2$, denote by $DF^i(r, 2s)$ the graph obtained from $DF(r, 2s)$ by merging all vertices in the same block of the partition of $V_i$ into a single vertex of degree $2(2r+1)$ and $3(2r+1)$, respectively.

\nt Suppose $2r+1=r_1r_2$ with $r_1, r_2\ge 3$. Let $V_3$ be the set of $r_1r_2$ vertices of degree $3s$ with induced vertex label $(21k+12)s$. Partition $V_3$ into $r_1\ge 3$ blocks of size $r_2\ge 3$. Clearly, none of these vertices have common neighbors. Denote by $DF^3(r_1, r, 2s)$ the graph obtained from $DF(r, 2s)$ by merging all vertices in the same block of the partition of $V_3$ into a single vertex of degree $3sr_2$.

\begin{theorem}\label{thm-DFr(2s)}  Let $s\ge 1$ be odd and $r\ge 1$. Suppose $(2r+1)s=2k+1$.
\begin{enumerate}[1.]
\item $\chi_{la}(DF^1(r, 2s))= 3$ if $k\not\equiv 2\pmod 4$;
\item $\chi_{la}(DF^2(r, 2s))= 3$;
\item $\chi_{la}(DF^3(r_1, r, 2s))= 3$, here $2r+1=r_1r_2$ with $r_1, r_2\ge 3$.
\end{enumerate}
\end{theorem}

\begin{proof}  Consider $DF^1(r, 2s)$. The $2s$ vertices of degree $2(2r+1)$ have induced vertex label $(2r+1)(10k+6)$, the degree 3 and $3s$ vertices are with induced vertex labels $9k+6$ and $(21k+12)s = (21k+12)(2k+1)/(2r+1)$. Clearly $9k+6$ is less than both $(2r+1)(10k+6)$ and $(21k+12)s$. Suppose $(2r+1)(10k+6)=(21k+12)s$. By assumption $(2r+1)s=2k+1$, we have $s$ is odd and $(2k+1)(10k+6)=(21k+12)s^2$. Thus, $(2k+1)(2k+2)\equiv k\pmod 4$.  This implies that $k\equiv 2\pmod 4$. So if $k\not\equiv 2\pmod 4$, then all vertex labels are distinct. We have $\chi_{la}(DF^1(r, 2s))\le 3$.

\ms\nt Consider $DF^2(r, 2s)$. The $s$ vertices of degree $3(2r+1)$ have induced vertex label $(2r+1)(9k+6)$, the degree 2 and $3s$ vertices are with induced vertex labels $10k+6$ and $(21k+12)s=(21k+12)(2k+1)/(2r+1)$ respectively. Since $2r+1\ge 3$, $10k+6$ is less than both $(2r+1)(9k+6)$ and $(21k+12)s$. Suppose $(2r+1)(9k+6)=(21k+12)s$. Similarly, we have $s$ is odd and $(2k+1)(9k+6)=(21k+12)s^2$. Thus $(2k+1)(k+2)\equiv k\pmod 4$. This implies that $k^2\equiv -1\pmod 4$ which is impossible.
 Since all the labels are distinct, $\chi_{la}(DF^2(r, 2s))\le 3$.

\ms\nt Consider $DF^3(r_1, r, 2s)$. The $r_1$ vertices of degree $3sr_2$ have induced vertex label $r_2(21k+12)s$, the degree 2 and 3 vertices are with induced vertex labels $10k+6$ and $9k+6$ respectively. Since all the labels are distinct, $\chi_{la}(DF^3(r_1, r, 2s))\le 3$.

\ms\nt Since the chromatic number of each considering graph is 3, the theorem holds.
\end{proof}

\section{New construction for $m=1$}

In this section, we give another way to construct a $5 \times (2k+1)$ matrix. We shall make use of the following table to obtain various families of graphs with triangles (and $C_6$ induced cycles) only. \vskip-8mm
\begin{table}[H]
\fontsize{9}{12}\selectfont
\begin{align*}&\begin{tabu}{|c|[1pt]c|c|c|c|c|c|c|[1pt]c|[1pt]l}\hline
i & 1 & 2 & 3 & \cdots  & k-2 & k-1 & k & k+1 &\\\tabucline[1pt]{-}
R_1 & 1 & 3 & 5 & \cdots  & 2k-5 & 2k-3 & 2k-1 & 2k+1 & \\\hline
R_2 & 4k+2 & 4k+1 & 4k & \cdots  & 3k+5 & 3k+4 & 3k+3 & 3k+2 & \\\hline
R_3  & 5k+3 & 5k+2 & 5k+1 & \cdots  & 4k+6 & 4k+5 & 4k+4 & 4k+3 &\\\hline
R_4 & 8k+4 & 8k+3 & 8k+2 & \cdots  & 7k+7 & 7k+6 & 7k+5 & 7k+4& \\\hline
R_5 & 8k+5 & 8k+7 & 8k+9 & \cdots & 10k-1  & 10k+1 & 10k+3 & 10k+5 &\\\hline
\end{tabu}\\
&\qquad \qquad
\begin{tabu}{|c|[1pt]c|c|c|c|c|c|c|}\hline
i &  k+2 & k+3 & k+4 & \cdots & 2k-1 & 2k & 2k+1 \\\tabucline[1pt]{-}
R_1 & 2 & 4 & 6 & \cdots & 2k-4 & 2k-2 & 2k \\\hline
R_2 & 3k+1 & 3k & 3k-1 & \cdots & 2k+4 & 2k+3 & 2k+2 \\\hline
R_3 & 6k+3 & 6k+2 & 6k+1 & \cdots & 5k+6 & 5k+5 & 5k+4 \\\hline
R_4 & 7k+3 & 7k+2 &  7k+1 & \cdots & 6k+6 & 6k+5 & 6k+4 \\\hline
R_5 & 8k+6 & 8k+8 & 8k+10 & \cdots & 10k & 10k+2 & 10k+4 \\\hline
\end{tabu}\end{align*}
\vskip-5mm\caption{}\label{table-Pt}
\end{table}

\normalsize
\nt We shall now describe a way to trace two sequences of numbers obtained from Table~\ref{table-Pt}.

\ms\nt Suppose $k$ is even. The first two terms of sequence $\mathcal S_1$ correspond to the $(R_2,k+1)$- and $(R_1,k+1)$-entry, followed by $4k$ terms correspond to the entries of the following consecutive $k/2$ segments
\begin{align*} (R_5,i), & (R_4,i), (R_2,2k+2-i), (R_1,2k+2-i),\\ & (R_5,k+1+i), (R_4,k+1+i), (R_2,k+1-i), (R_1,k+1-i), \quad 1\le i\le \frac{k}{2}.\end{align*}

\nt Similarly, the first two terms of sequence $\mathcal S_2$ correspond to the $(R_4,k+1)$- and $(R_5,k+1)$-entry, followed by $4k$ terms correspond to the entries of the following consecutive $k/2$ segments
\begin{align*} (R_1,j), & (R_2,j), (R_4,2k+2-j), (R_5,2k+2-j),\\ & (R_1,k+1+j), (R_2,k+1+j), (R_4,k+1-j), (R_5,k+1-j),\quad 1\le j\le \frac{k}{2}.\end{align*}
\nt Note that the $(R_1, \frac{k}{2}+1)$-entry is $k+1$ and the $(R_5, \frac{k}{2}+1)$-entry is $9k+5$.

\ms\nt We now have the following observations.
\begin{enumerate}[(A)]
  \item The sum of the first (respectively, the last) term of $\mathcal S_1$ and $\mathcal S_2$ is $10k+6$.
  \item In both $\mathcal S_1$ and $\mathcal S_2$, the sum of the $2r$-th and $(2r+1)$-st terms is also $10k+6$ for $r\in [1, 2k]$.
  \item For $1\le j\le 2k$, the $(2j-1)$-st and $2j$-th terms are in the same column of the table. Together with the third row entry of the same column, the sum is either $9k+6$ (if they belong to the first 3 rows) or else is $21k+12$ (if they belong to the last 3 rows).
\end{enumerate}

\ms\nt Suppose $k$ is odd. We also can obtain two sequences of numbers from the table above  similarly having the same properties. The first two terms of $\mathcal S_1$ are $(R_2,k+1)$, $(R_1, k+1)$ followed by $4(k-1)$ terms correspond to the entries of the following consecutive $(k-1)/2$ segments
\begin{align*} (R_5,i), & (R_4,i), (R_2,2k+2-i), (R_1,2k+2-i), (R_5,k+1+i),\\&  (R_4,k+1+i), (R_2,k+1-i), (R_1,k+1-i),\quad 1\le i\le \frac{k-1}{2},\end{align*}
and then the sequence ends with last four terms $(R_5,\frac{k+1}{2})$, $(R_4,\frac{k+1}{2})$, $(R_2,\frac{3k+3}{2})$, $(R_1,\frac{3k+3}{2})$.

\ms\nt Similarly, the first two terms of $\mathcal S_2$ are $(R_4,k+1)$, $(R_5,k+1)$ followed by $4(k-1)$ terms correspond to the entries of the following consecutive $(k-1)/2$ segments
\begin{align*} (R_1,j), & (R_2,j), (R_4,2k+2-j), (R_5,2k+2-j),\\ & (R_1,k+1+j), (R_2,k+1+j), (R_4,k+1-j), (R_5,k+1-j), \quad 1\le j\le \frac{k-1}{2},\end{align*}
and then the sequence ends with last four terms $(R_1,\frac{k+1}{2})$, $(R_2,\frac{k+1}{2})$, $(R_4,\frac{3k+3}{2})$, $(R_5,\frac{3k+3}{2})$.

\ms\nt For $n\ge 1$, we are ready to define a {\it peanut graph} that consists of two $3$-cycles and $n$ $6$-cycles, denoted $Pt(n)$, with vertex set \[\{x, y\}\cup\{u_i, v_i\;|\; 1\le i\le 2n+1\}\] and edge set \[\{xu_1, xv_1, yu_{2n+1}, yv_{2n+1}\}\cup\{u_iu_{i+1},v_iv_{i+1}\;|\; 1\le i\le 2n\}\cup\{u_{2j-1}v_{2j-1}\;|\; 1\le j\le n+1\}.\] Note that $Pt(n)$ has order $4n+4$ and size $5n+5$. Moreover, vertices $x,y$ and $u_{2i}, v_{2i}$ are of degree 2 for $i\in [1,n]$ and vertices $u_{2j-1}v_{2j-1}$ are of degree 3 for $j\in [1,n+1]$.

\begin{theorem}\label{thm-Pt(n)even} For even $n\ge 2$, $\chi_{la}(Pt(n)) = 3$. \end{theorem}

\begin{proof} Let $n=2k\ge 2$. Since $Pt(2k)$ contains $3$-cycles, $\chi_{la}(Pt(2k))\ge \chi(Pt(2k))=3$. It suffices to give a local antimagic 3-coloring of $Pt(2k)$. For convenience, we denote $y$ by $u_{4k+2}=v_{4k+2}$.

\ms\nt Note that $Pt(2k)$ has two induced paths of order $4k+3$, namely $\mathcal P_1 = xu_1u_2u_3\cdots u_{4k+1}u_{4k+2}$ and $\mathcal P_2 = xv_1v_2v_3\cdots v_{4k+1}v_{4k+2}$.

\ms\nt Suppose $k\ge 2$ is even. We shall define a bijection $f: E(Pt(2k))\to [1,10k+5]$ such that the consecutive edges of $\mathcal P_1$ (respectively, $\mathcal P_2$) are labeled by the terms of $\mathcal S_1$ (respectively $\mathcal S_2$) correspondingly.

\ms\nt Functionally, for $1\le i\le k/2$, $f_1: E(\mathcal P_1)\to \{2k+1,3k+2\}\cup [k+1,2k] \cup [2k+2,5k/2+1] \cup [3k+3,7k/2+2] \cup [13k/2+4,7k+3] \cup [8k+5,9k+4]$ such that
\begin{multicols}{2}
\begin{enumerate}[(1)]
  \item $f_1(xu_1) = 3k+2, f_1(u_1u_2) = 2k+1$; 
  \item $f_1(u_{8i-6}u_{8i-5}) = 8k+3+2i$; 
  \item $f_1(u_{8i-5}u_{8i-4}) = 8k+5-i$; 
  \item $f_1(u_{8i-4}u_{8i-3}) = 2k+1+i$; 
  \item $f_1(u_{8i-3}u_{8i-2}) = 2k+2-2i$; 
  \item $f_1(u_{8i-2}u_{8i-1}) = 8k+4+2i$; 
  \item $f_1(u_{8i-1}u_{8i}) = 7k+4-i$; 
  \item $f_1(u_{8i}u_{8i+1}) = 3k+2+i$;  
  \item $f_1(u_{8i+1}u_{8i+2}) = 2k+1-2i$. 
\end{enumerate}
\end{multicols}
\nt Moreover, for $1\le i\le k/2$, $f_2 : E(\mathcal P_2)\to \{7k+4,10k+5\}\cup [1,k] \cup [5k/2+2,3k+1] \cup [7k/2+3, 4k+2]\cup [6k+4, 7k+3]\cup [7k+5,8k+4] \cup [9k+5,10k+4]$ such that
\begin{multicols}{2}
\begin{enumerate}[(1)]
  \item $f_2(xv_1) = 7k+4, f_2(v_1v_2) = 10k+5$;
 \item $f_2(v_{8i-6}v_{8i-5}) = 2i-1$; 
  \item $f_2(v_{8i-5}v_{8i-4}) = 4k+3-i$; 
  \item $f_2(v_{8i-4}v_{8i-3}) = 6k+3+i$; 
  \item $f_2(v_{8i-3}v_{8i-2}) = 10k+6-2i$; 
  \item $f_2(v_{8i-2}v_{8i-1}) = 2i$; 
  \item $f_2(v_{8i-1}v_{8i}) = 3k+2-i$; 
  \item $f_2(v_{8i}v_{8i+1}) = 7k+4+i$;  
  \item $f_2(v_{8i+1}v_{8i+2}) = 10k+5-2i$. 
\end{enumerate}
\end{multicols}

\ms\nt Thus, we have used the entries of rows 1, 2, 4, 5 of Table~\ref{table-Pt}, i.e., integers in $[1,4k+2]\cup [6k+4,10k+5]$. We note that edges
\begin{enumerate}[(i)]
  \item $xu_1, u_1u_2, xv_1, v_1v_2$ are labeled by entries in the $(k+1)$-st column,
  \item $u_{2j-2}u_{2j-1}, u_{2j-1}u_{2j}, v_{2j-2}v_{2j-1}, v_{2j-1}v_{2j}$, $2\le j\le n$, are labeled by edges in the same column, namely, column $1, 2k+1, k+2, k; 2, 2k, k+3, k-1; 3, 2k-1, k+4, k-2; \ldots; k/2-1, 3k/2+3, 3k/2, k/2+2; k/2, 3k/2+2, 3k/2+1$.
  \item $yu_{2n+1}, u_{2n}u_{2n+1}, yv_{2n+1}, v_{2n}v_{2n+1}$ are labeled by entries in the $(k/2+1)$-st column.
\end{enumerate}

\nt In each of (i) to (iii) above, every group of 4 edges are adjacent to the edge $u_{2j-1}v_{2j-1}$ respectively for $j\in [1,n+1]$. Thus, we label the edge $u_{2j-1}v_{2j-1}$ by the row 3 entry of the corresponding column respectively. Functionally, we have $f_3 : \{u_{2j-1}v_{2j-1}\;|1\le j\le 2k+1\}\to [4k+3,6k+3]$ such that for $1\le j\le k/2$,
\begin{multicols}{2}
\begin{enumerate}[(1)]
  \item $f_3(u_1v_1) = 4k+3$;
  \item $f_3(u_{8j-5}v_{8j-5}) = 5k+4-j$; 
  \item $f_3(u_{8j-3}v_{8j-3}) = 5k+3+j$; 
  \item $f_3(u_{8j-1}v_{8j-1}) = 6k+4-j$; 
  \item $f_3(u_{8j+1}v_{8j+1}) = 4k+3+j$. 
  \item[]
\end{enumerate}
\end{multicols}

\nt Thus, by Observation $(A)$ to $(C)$ above, combining $f_1$, $f_2$ and $f_3$ we obtain a local antimagic 3-coloring, say $f$, of $Pt(2k)$ such that every degree 2 (respectively, degree 3) vertex of $Pt(2k)$ has induced label $10k+6$ (respectively,  $9k+6$ and $21k+12$ alternately if in $\mathcal P_1$; and $21k+12$ and $9k+6$ alternately if in $\mathcal P_2$). Thus, $f$ is a local antimagic 3-coloring of $Pt(2k)$.

\ms\nt Suppose $k\ge 1$ is odd. The corresponding functions $f_1, f_2, f_3$ are given as follows for $1\le i,j\le (k+1)/2$.
\begin{multicols}{2}
\begin{enumerate}[(1)]
  \item $f_1(xu_1) = 3k+2, f_1(u_1u_2) = 2k+1$; 
  \item $f_1(u_{8i-6}u_{8i-5}) = 8k+3+2i$; 
  \item $f_1(u_{8i-5}u_{8i-4}) = 8k+5-i$; 
  \item $f_1(u_{8i-4}u_{8i-3}) = 2k+1+i$; 
  \item $f_1(u_{8i-3}u_{8i-2}) = 2k+2-2i$; 
  \item $f_1(u_{8i-2}u_{8i-1}) = 8k+4+2i$, $i\ne (k+1)/2$; 
  \item $f_1(u_{8i-1}u_{8i}) = 7k+4-i$, $i\ne (k+1)/2$; 
  \item $f_1(u_{8i}u_{8i+1}) = 3k+2+i$, $i\ne (k+1)/2$;  
  \item $f_1(u_{8i+1}u_{8i+2}) = 2k+1-2i$,  $i\ne (k+1)/2$. 
  \item[]
\end{enumerate}
\end{multicols}
\begin{multicols}{2}
\begin{enumerate}[(1)]\addtocounter{enumi}{9}
  \item $f_2(xv_1) = 7k+4, f_2(v_1v_2) = 10k+5$;
 \item $f_2(v_{8i-6}v_{8i-5}) = 2i-1$; 
  \item $f_2(v_{8i-5}v_{8i-4}) = 4k+3-i$; 
  \item $f_2(v_{8i-4}v_{8i-3}) = 6k+3+i$; 
  \item $f_2(v_{8i-3}v_{8i-2}) = 10k+6-2i$; 
  \item $f_2(v_{8i-2}v_{8i-1}) = 2i$, $i\ne (k+1)/2$; 
  \item $f_2(v_{8i-1}v_{8i}) = 3k+2-i$, $i\ne (k+1)/2$; 
  \item $f_2(v_{8i}v_{8i+1}) = 7k+4+i$, $i\ne (k+1)/2$;  
  \item $f_2(v_{8i+1}v_{8i+2}) = 10k+5-2i$, $i\ne (k+1)/2$. 
  \item[]
  \end{enumerate}
  \end{multicols}
  \begin{multicols}{2}
\begin{enumerate}[(1)]\addtocounter{enumi}{18}
  \item $f_3(u_1v_1) = 4k+3$;
  \item $f_3(u_{8j-5}v_{8j-5}) = 5k+4-j$; 
  \item $f_3(u_{8j-3}v_{8j-3}) = 5k+3+j$; 
  \item $f_3(u_{8j-1}v_{8j-1}) = 6k+4-j$, $j\ne (k+1)/2$; 
  \item $f_3(u_{8j+1}v_{8j+1}) = 4k+3+j$, $j\ne (k+1)/2$. 
  \item[]
\end{enumerate}
\end{multicols}

\nt By a similar argument, we can also conclude that $Pt(2k)$ is local antimagic 3-colorable. This completes the proof.
\end{proof}

\begin{example} We give the graphs $Pt(4)$ and $Pt(10)$, in Figures~\ref{fig:Pt4} and~\ref{fig:Pt10} respectively, with the labelings given by the tables below.

\begin{multicols}{2}
\nt For $k=2$, the table is
\[\fontsize{8}{11}\selectfont
\begin{tabu}{|c|[1pt]c|c|[1pt]c|[1pt]c|c|}\hline
R_1 & 1 & 3 & 5 & 2 & 4  \\\hline
R_2 & 10 &9 & 8 & 7 & 6  \\\hline
R_3  & 13 & 12 & 11 & 15 & 14 \\\hline
R_4 & 20 & 19 & 18 & 17 & 16 \\\hline
R_5 & 21 & 23 & 25  & 22 & 24 \\\hline
\end{tabu}
\]
\nt and $\mathcal S_1 = 8, 5, 21, 20, 6, 4, 22, 17, 9, 3$,\\ \hspace*{6.5mm} $\mathcal S_2 = 18, 25, 1, 10, 16, 24, 2, 7, 19, 23$. 

\begin{figure}[H]
\begin{center}\vskip1cm
\epsfig{file=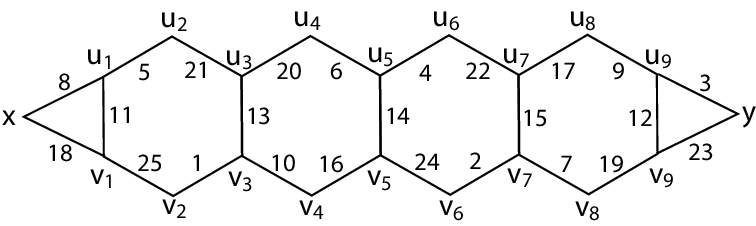, width=7.6cm}
\caption{Graph $Pt(4)$.}\label{fig:Pt4}
\end{center}
\end{figure}
\end{multicols}

\ms\nt For $k=5$, the table is
\[\fontsize{8}{11}\selectfont
\begin{tabu}{|c|[1pt]c|c|c|c|c|[1pt]c|[1pt]c|c|c|c|c|}\hline
R_1 & 1 & 3 & 5 & 7 & 9 & 11 & 2 & 4 & 6 & 8 & 10 \\\hline
R_2 & 22 & 21 & 20 & 19 & 18 & 17 & 16 & 15 & 14 & 13 & 12 \\\hline
R_3 & 28 & 27 & 26 & 25 & 24 & 23 & 33 & 32 & 31 & 30 & 29  \\\hline
R_4 & 44 & 43 & 42 & 41 & 40 & 39 & 38 & 37 & 36 & 35 & 34 \\\hline
R_5 & 45 & 47 & 49 & 51 & 53 & 55 & 46 & 48 & 50 & 52 & 54 \\\hline
\end{tabu}\]
and\\
$\begin{aligned}\mathcal S_1 & = 17, 11, 45, 44, 12, 10, 46, 38, 18, 9, 47, 43, 13, 8, 48, 37, 19, 7, 49, 42, 14,6,\\
\mathcal S_2 & = 39, 55, 1, 22, 34, 54, 2, 16, 40, 53, 3, 21, 35, 52, 4, 15, 41, 51, 5, 20, 36, 50.
\end{aligned}$\\
The sequence of labels for edges $u_{2j-1}v_{2j-1}$, $1\le j\le 11$, is 23, 28, 29, 33, 24, 27, 30, 32, 25, 26, 31.\\
\vskip-0.5cm
\begin{figure}[H]
\centerline{\epsfig{file=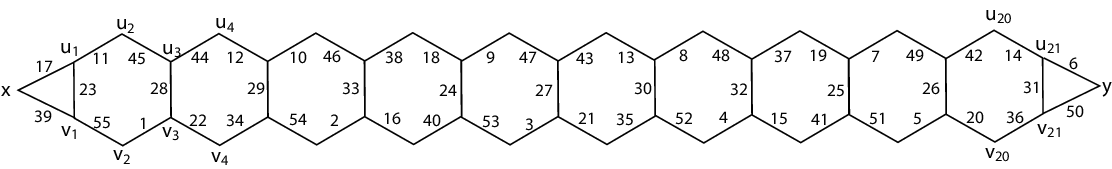, width=14cm}}
\caption{Graph $Pt(10)$.}\label{fig:Pt10}
\end{figure}
\vskip-0.5cm\rsq
 \end{example}

\nt Let $TB(n)$ be the {\it triangular bracelet} obtained from $Pt(n)$ by merging vertices $x$ and $y$ to get vertex $z_0$, and vertices $u_{2i}$ and $v_{2i}$ to get vertex $z_{2i}$ for $1\le i\le n$. Note that $TB(n)$ has $2n+2$ triangles having order $3n+3$ and size $5n+5$. There are $n+1$ vertices of degree 4 and $2n+2$ vertices of degree 3.

\begin{theorem}\label{thm-TBn} For even $n\ge 2$, $\chi_{la}(TB(n))=3$. \end{theorem}

\begin{proof} Let $n=2k\ge 2$. Clearly, $\chi_{la}(TB(2k))\ge 3$. By keeping the same edge labeling for $Pt(2k)$, we immediately have a local antimagic 3-coloring $f: E(TB(2k))\to [1,10k+5]$ such that every two adjacent degree 3 vertices has induced label $9k+6$ or $21k+12$ while every degree 4 vertex has induced label $20k+12$. This completes the proof. \end{proof}

\begin{example}\label{ex-TB10}
By merging $x$ with $y$, and $u_{2j}$ with $v_{2j}$, $1\le j\le 10$, from $Pt(10)$, we have a local antimagic $3$-coloring for $TB(10)$ shown as follows:
\begin{figure}[H]
\centerline{\epsfig{file=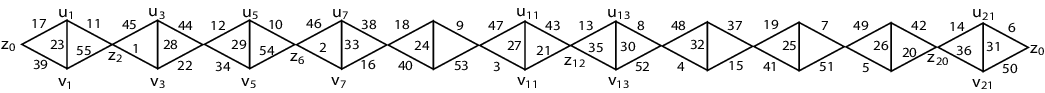, width=14cm}}
\caption{Graph $TB(10)$.}\label{fig:TB10}
\end{figure}
\end{example}

\nt By an argument similar to the proof of Theorem~\ref{thm-FBrs12}, we also have the following theorem.

\begin{theorem}\label{thm-PtTBgeneral} Consider even $n\ge 2$.
\begin{enumerate}[(a)]
  \item  Suppose $V_1$ (resp. $V_2$) is the set of $n+1$ degree $3$ vertices of $Pt(n)$ with induced vertex label $9k+6$ (resp. $21k+12$). Partition $V_1$ (resp. $V_2$) into $r_1$ block(s)  (resp. $r_2$ block(s)) each of size $s_1=(n+1)/r_1$ (resp. $s_2=(n+1)/r_2$) such that $r_i\ge 1, s_i\ge 3$ are both odd. If $Pt^i(r_i,n)$ is the graph obtained from $Pt(n)$ by merging all the degree $3$ vertices in the same block of the partition of $V_i$, then $\chi_{la}(Pt^i(r_i,n))=3$, $i=1,2$. 
  \item Suppose $V_3$ is the set of $2n+2$ degree $2$ vertices of $Pt(n)$. Partition $V_3$ into $r\ge 2$ blocks of equal size $2\le s = (2n+2)/r\le n+1$ such that all vertices in the same block have no common neighbors. If $Pt^3(r,n)$ is the graph obtained from $Pt(n)$ by merging all vertices in the same block of the partition of $V_3$, then $\chi_{la}(Pt^3(r,n))=3$. 
\item Suppose $U_1$ (resp. $U_2$) is the set of $n+1$ degree $3$ vertices of $TB(n)$ with induced vertex label $9k+6$ (resp. $21k+12$). Partition $V_1$ (resp. $V_2$) into $r_1$ block(s)  (resp. $r_2$ block(s)) each of size $s_1=(n+1)/r_1$ (resp. $s_2=(n+1)/r_2$) such that $r_i, s_i\ge 3$ are both odd are all vertices in the same block have no common neighbors.  If $TB^i(r_i, n)$ is the graph obtained from $TB(n)$ by merging all the degree $3$ vertices in the same block of the partition of $U_i$, then $\chi_{la}(TB^i(r_i,n))=3$. 
  \item Suppose $U_3$ is the set of $n+1$ degree $4$ vertices of $TB(n)$, $n\ge 8$. Partition $U_3$ into $r$ blocks of odd size $s=(n+1)/r$ such that $r,s\ge 3$ are both odd and all vertices in the same  block have no common neighbors. Let $TB^3(r,n)$ be the graph obtained by merging all vertices in the same block, then $\chi_{la}(TB^3(r,n))=3$. 
\end{enumerate}
\end{theorem}

\begin{proof} Let $n=2k$. Keep the labeling of $Pt(n)$ as defined in the proof of Theorem~\ref{thm-Pt(n)even}. Now, every degree 2 (respectively, degree 3) vertex of $Pt(2k)$ has induced label $10k+6$ (respectively,  $9k+6$ and $21k+12$. Also keep the labeling of $TB(n)$ as defined in the proof of Theorem~\ref{thm-TBn}. Now, every degree 4 vertex has induced label $20k+12$ while every two adjacent degree 3 vertices have labels $9k+6$ or $21k+12$.

\begin{enumerate}[(a)]
\item  For $Pt^1(r_1,n)$, the induced vertex labels are now $10k+6$, $s_1(9k+6)$, $21k+12$. For $Pt^2(r_2,n)$, the induced vertex labels are now $10k+6$, $9k+6$, $s_2(21k+12)$. Thus, $Pt^i(r_i,n)$, $i=1,2,$ are local antimagic 3-colorable since $s_1\ge 3$.

\item  For $Pt^3(r,n)$, the induced vertex labels are now $s(10k+6)$, $9k+6$, $21k+12$ that are all distinct since $s\ge 2$. Thus, $Pt^3(r,n)$ is local antimagic 3-colorable.

\item  For $TB^1(r_1,n)$, the induced vertex labels are now $s_1(9k+6)$, $21k+12$, $20k+12$. For $TB^2(r_2,n)$, the induced vertex labels are now $s_2(21k+12)$, $9k+6$, $20k+12$. Thus, $TB^i(r_i,n)$, $i=1,2,$ are local antimagic 3-colorable since $s_1\ge 3$.

\item  For $TB^3(r,n)$, the induced vertex labels are now $s(20k+12)$, $9k+6$, $21k+12$ that are all distinct since $s\ge 3$. Thus, $TB^3(r,n)$ is local antimagic 3-colorable.
\end{enumerate}
\ms\nt Since all the graphs obtained have $\chi = 3$, the theorem holds.\end{proof}

\nt Note that the graphs in Theorem~\ref{thm-PtTBgeneral} above may not be unique, and the graph in (b) is a generalization of $TB(n)$. We can now give a theorem on disjoint union of triangular bracelets.

\ms\nt Consider the graph $TB(s)+ TB(4i_1-2) + TB(4i_2-2) + \cdots + TB(4i_r-2)$, where $1\le i_1 < i_2 < \cdots < i_r$. The order of this graph is $3[s+1+\sum\limits_{a=1}^r (4i_a-1)]=3[s+\sum\limits_{a=1}^r (4i_a-1)]+3$. We denote this graph by $G(n)$, where $n=s+\sum\limits_{a=1}^r (4i_a-1)$. So the order of $G(n)$ is $3(n+1)$.

\begin{theorem}\label{thm-Gn} Let $G(n)$ be the graph defined above satisfying the following conditions:
\begin{enumerate}[(a)]
  \item $8i_{a+1} > 16i_a-2$ for $1\le a\le r-1$.
  \item $n\ge 8i_r-2$.
\end{enumerate}
Then $\chi_{la}(G(n))=3$. \end{theorem}

\begin{proof} Clearly, $\chi_{la}(G(n))\ge \chi(G(n)) = 3$. Suffice to show that $G(n)$ admits a local antimagic 3-coloring. Suppose $n=2k\ge 2$. Begin with $TB(n)$ and its local antimagic 3-coloring $f$, which is the combined labeling of $f_1$, $f_2$ and $f_3$, as defined in the proof of Theorem~\ref{thm-TBn}. We note that the edges given by the cycle
\[\mathcal C_1=z_0u_1z_2u_3z_4\cdots u_{2n-1}z_{2n}u_{2n+1}z_0\] (respectively, $\mathcal C_2=z_0v_1z_2v_3z_4\cdots v_{2n-1}z_{2n}v_{2n+1}z_0$) has consecutive edge labels given by $g_1 = f_1$ (respectively, $g_2=f_2$) by replacing vertices $x,y$ by $z_0$, and vertices $u_{2i}, v_{2i}$ by $z_{2i}$, $1\le i\le n$. For convenience, let $z_{2(n+1)}=z_0$.

\nt Suppose $k$ is even. For clarity, we list $g_1,g_2$ as follows with $1\le i\le k/2$:

\begin{multicols}{2}
\begin{enumerate}[(1)]
  \item $g_1(z_0u_1) = 3k+2$;
  \item $g_1(u_1z_2) = 2k+1$; 
  \item $g_1(z_{8i-6}u_{8i-5}) = 8k+3+2i$; 
  \item $g_1(u_{8i-5}z_{8i-4}) = 8k+5-i$;
  \item $g_1(z_{8i-4}u_{8i-3}) = 2k+1+i$;
  \item $g_1(u_{8i-3}z_{8i-2}) = 2k+2-2i$;
  \item $g_1(z_{8i-2}u_{8i-1}) = 8k+4+2i$;
  \item $g_1(u_{8i-1}z_{8i}) = 7k+4-i$;
  \item $g_1(z_{8i}u_{8i+1}) = 3k+2+i$;
  \item $g_1(u_{8i+1}z_{8i+2}) = 2k+1-2i$;
\end{enumerate}
\begin{enumerate}[(1)]
  \item $g_2(z_0v_1) = 7k+4$;
  \item $g_2(v_1z_2) = 10k+5$;
  \item $g_2(z_{8i-6}v_{8i-5}) = 2i-1$; 
  \item $g_2(v_{8i-5}z_{8i-4}) = 4k+3-i$; 
  \item $g_2(z_{8i-4}v_{8i-3}) = 6k+3+i$; 
  \item $g_2(v_{8i-3}z_{8i-2}) = 10k+6-2i$; 
  \item $g_2(z_{8i-2}v_{8i-1}) = 2i$; 
  \item $g_2(v_{8i-1}z_{8i}) = 3k+2-i$; 
  \item $g_2(z_{8i}v_{8i+1}) = 7k+4+i$;  
  \item $g_2(v_{8i+1}z_{8i+2}) = 10k+5-2i$. 
\end{enumerate}
\end{multicols}

\nt We now have
\begin{multicols}{2}\fontsize{10}{10}\selectfont
\begin{enumerate}[(1)]
\item $g_1(z_0u_1) + g_2(z_0v_1) = 10k+6$;
\item $g_1(u_1z_2) + g_2(v_1z_2) =  12k+6$;
\item $g_1(z_{8i-6}u_{8i-5}) +  g_2(z_{8i-6}v_{8i-5}) = 8k+2+4i$;
\item $g_1(u_{8i-5}z_{8i-4}) + g_2(v_{8i-5}z_{8i-4}) =  12k+8-2i$;
\item $g_1(z_{8i-4}u_{8i-3}) + g_2(z_{8i-4}v_{8i-3})= 8k+4+2i$;
\item $g_1(u_{8i-3}z_{8i-2}) + g_2(v_{8i-3}z_{8i-2}) = 12k+8-4i$;
\item $g_1(z_{8i-2}u_{8i-1}) + g_2(z_{8i-2}v_{8i-1})= 8k+4+4i$;
\item $g_1(u_{8i-1}z_{8i}) + g_2(v_{8i-1}z_{8i})=10k+6-2i$;
\item $g_1(z_{8i}u_{8i+1}) + g_2(z_{8i}v_{8i+1}) = 10k+6+2i$;
\item $g_1(u_{8i+1}z_{8i+2}) + g_2(v_{8i+1}z_{8i+2})= 12k+6-4i$.
\end{enumerate}\normalsize
\end{multicols}

\nt From (4) and (5) we have
\begin{eqnarray*}
g_1(u_{8j-5}z_{8j-4}) + g_2(v_{8j-5}z_{8j-4}) & = & 12k+8-2j;\\
g_1(z_{8j-4}u_{8j-3}) + g_2(z_{8j-4}v_{8j-3}) & = & 8k+4+2j,
\end{eqnarray*}
while from (6) and (7) we have
\begin{eqnarray*}
g_1(u_{8i-3}z_{8i-2}) + g_2(v_{8i-3}z_{8i-2}) &=& 12k+8-4i;\\
g_1(z_{8i-2}u_{8i-1}) + g_2(z_{8i-2}v_{8i-1}) &=& 8k+4+4i.
\end{eqnarray*}
\nt Thus, $12k+8-2j = 12k+8-4i$ and $8k+4+2j=8k+4+4i$ if and only if $j=2i$ for $i\ge 1$. This implies that
\begin{alignat*}{3}
g_1(u_{8i-3}z_{8i-2}) + g_2(v_{8i-3}z_{8i-2})  & = & g_1(u_{16i-5}z_{16i-4}) + g_2(v_{16i-5}z_{16i-4})& = & 12k+8-4i,\\
g_1(z_{8i-2}u_{8i-1}) + g_2(z_{8i-2}v_{8i-1})  & = & g_1(z_{16i-4}u_{16i-3}) + g_2(z_{16i-4}v_{16i-3}) & = & 8k+4+4i.
\end{alignat*}

\nt Consider the integers $i_1, i_2, \ldots, i_r$ $(r\ge 1)$, where $1\le i_1 < i_2 < \cdots < i_r$ such that
\begin{enumerate}[(a)]
  \item $8i_{a+1} > 16i_a-2$ (equivalently, $8i_{a+1} -2 > 16i_a-4$) for $1\le a\le r-1$.
  \item $n \ge 8i_r-2$.
\end{enumerate}

\nt For $1\le a\le r$, we can now split vertex $z_{8i_a-2}$ (respectively, $z_{16i_a-4}$) into vertices $z^1_{8i_a-2}$ and $z^2_{8i_a-2}$ (respectively, $z^1_{16i_a-4}$ and $z^2_{16i_a-4}$) such that
\begin{enumerate}[(i)]
  \item $z^1_{8i_a-2}$ is incident to edges with labels $2k+2-2i_a$ and $10k+6-2i_a$, that is, $z^1_{8i_a-2}$ is adjacent to $u_{8i_a-3}$ and $v_{8i_a-3}$,
  \item $z^2_{8i_a-2}$ is incident to edges with labels $8k+4+2i_a$ and $2i_a$, that is,  $z^2_{8i_a-2}$ is adjacent to $u_{8i_a-1}$ and $v_{8i_a-1}$,
  \item $z^1_{16i_a-4}$ is incident to edges with labels $8k+5-2i_a$ and $4k+3-2i_a$, that is,  $z^1_{16i_a-4}$ is adjacent to $u_{16i_a-5}$ and $v_{16i_a-5}$,
  \item $z^2_{16i_a-4}$ is incident to edges with labels $2k+1+2i_a$ and $6k+3+2i_a$, that is, $z^2_{16i_a-4}$ is adjacent to $u_{16i_a-3}$ and $v_{16i_a-3}$.
\end{enumerate}
\nt Note that, condition~(a) makes sure that the sequence of indices
 $8i_1-2$, $16i_1-4$, $8i_2-2$, $16i_2-4$, $\dots$, $8i_r-2$, $16i_r-4$ is a strictly increasing sequence;  condition~(a) makes sure that $16i_r-4<2n+1$. Thus, $\mathcal C_1$ (also $\mathcal C_2$) is separated into $r+1$ cycles.

\ms\nt Finally, merge vertices $z^1_{8i_a-2}$ and $z^2_{16i_a-4}$ (respectively, $z^2_{8i_a-2}$ and $z^1_{16i_a-4}$) to get new vertices of degree 4 (still denote by $z_{8i_a-2})$ (respectively, $z_{16i_a-4}$) with induced vertex labels $20k+12$.
 Note that all other vertices of degree $4$ whose labels are $20k+12$ too and all vertices of degree $3$ whose labels unchange with labels $9k+6$ or $21k+12$, please see Theorem~\ref{thm-TBn}. Thus, $G$ is local antimagic 3-colorable.

\ms\nt Suppose $k$ is odd. By the same argument, we also get the same conclusion. This completes the proof. \end{proof}

\begin{example} Consider $TB(30)$. Under the condition~(b) in Theorem~\ref{thm-Gn}, $i_r\le 32/8=4$, $r\ge 1$.  The consecutive edge labels of $$\mathcal C_1= z_0u_1z_2u_3z_4\cdots u_{59}z_{60}u_{61}z_0 \mbox{ and } \mathcal C_2 = z_0v_1z_2v_3z_4\cdots v_{59}z_{60}v_{61}z_0$$ are:

\[\mathcal C_1: 47,31,125,124,32,{\bf30,126},108,48,29,127,{\bf123,33}, \underline{28,128},107,49,27,\]
\[129,122,34,\overline{26,130},106,50,25,131,\underline{121,35},{\it 24,132},105,51,23,133,120,36,22,134,104,52,21,\]
\[135,\overline{119,37},20,136,103,53,19,137,118,38,18,138,102,54,17,139,{\it 117,39},16, 47\]

\nt and

\[\mathcal C_2: 109,155,1,62,94,{\bf 154,2}, 46,110,153,3, {\bf 61,95}, \underline{152,4}, 45,111,151,\]
\[5,60,96,\overline{150,6},44,112,149,7,\underline{59,97},{\it 148,8},43,113,147,9,58,98,146,10,42,114,145,\]
\[11,\overline{57,99},144,12,41,115,143,13,56,100,142,14,40,116,141,15,{\it 55,101},140,109\]

\nt The corresponding edge labels of $u_{2j-1}v_{2j-1}$ for $1\le j\le 31$ are
\[63, 78,79,93,64, 77,80,92,65, 76,81,91,66, 75,82,90,67, \]\[74,83,89,68, 73,84,88,69, 72,85,87,70, 71, 86.\]

\begin{enumerate}[(1)]
  \item If $r=1$, then $i_1\in[1,4]$;
  \item If $r=2$, then $(i_1,i_2)\in\{(1,2),(1,3),(1,4),(2,4)\}$;
  \item If $r=3$, then $(i_1,i_2,i_3)=(1,2,4)$.
\end{enumerate}

\nt In (1), we can split vertices $z_6, z_{12}$ or $z_{14}, z_{28}$ or $z_{22}, z_{44}$ or $z_{30}, z_{60}$.

\ms\nt In (2), we can split vertices $z_6, z_{12}$ and $z_{14}, z_{28}$; or $z_6, z_{12}$ and $z_{22}, z_{44}$; or $z_6, z_{12}$ and $z_{30}, z_{60}$; or $z_{14}, z_{28}$ and $z_{30}, z_{60}$.

\ms\nt In (3), we can only split vertices $z_6, z_{12}$; $z_{14}, z_{28}$ and $z_{30}, z_{60}$. For this case, we obtain the graph $TB(5)+TB(2)+TB(6)+TB(14)$.

\ms\nt For easy illustration of the method in the proof of Theorem~\ref{thm-Gn}, we use a smaller $n$, say 10, and show at the next example. \rsq
\end{example}

\begin{example} Let $n=10$. Starting from $TB(10)$ with the labeling shown at Example~\ref{ex-TB10}. Under the condition (b) in Theorem~\ref{thm-Gn}, $r=1$ and $i_1=1$. Hence $s=7$. So we spit $z_6$ and $z_{12}$. According to the construction in the proof of Theorem~\ref{thm-Gn} we get the following local antimagic $3$-coloring for $G(10)=TB(7)+TB(2)$.
\begin{figure}[H]
\centerline{\epsfig{file=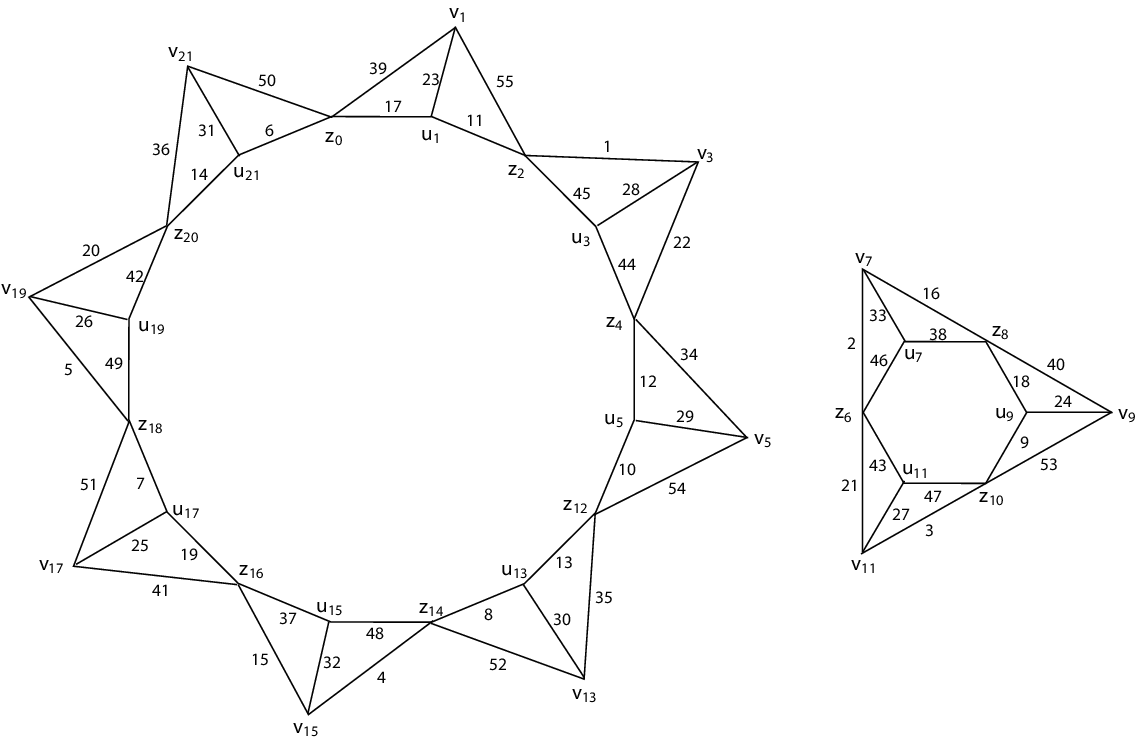, width=14cm}}
\caption{Graph $TB(7)+TB(2)$.}\label{fig:G10}
\end{figure}\rsq
\end{example}

\nt Recall that both $TB(n)$ and $G(n)$ above have $n+1$ degree 4 vertices. For even $n\ge 8$ such that $n+1 = rs$, $r,s\ge 3$, partition the $rs$ degree 4 vertex set into $r$ blocks of size $s$, in which vertices in the same block are mutually without common neighbors. Let {\it generalized bracelet}, denote $GB(n)$, be a graph obtained from $TB(n)$ or $G(n)$ by merging all the vertices in the same block. Note that $GB(n)$ has no multiple edges, may not be connected, and is not unique in general.

\begin{theorem} For even $n\ge 8$, $\chi_{la}(GB(n)) = 3$. \end{theorem}

\begin{proof} Consider the $TB(n)$ (or $G(n)$) with the local antimagic 3-coloring in Theorem~\ref{thm-TBn} (or Theorem~\ref{thm-Gn}) having degree 3 vertices induced label $9k+6$ or $21k+12$, and degree 4 vertex induced label $20k+12$. Clearly, $GB(n)$ also admits a local antimagic 3-coloring with degree 3 vertices having induced label $9k+6$ or $21k+12$, and there are $r$ degree $4s$ vertices with induced vertex label $s(20k+12)$. Thus, $\chi_{la}(GB(n))\le 3$. Since $\chi_{la}(GB(n))\ge \chi(GB(n))=3$, the theorem holds.
\end{proof}

\section{$m=3$}

In this section, we shall make use of the following $11\times (2k+1)$ matrix (with integers in $[1,22k+11]$ bijectively for $k\ge 1$) to get our first result.
\[\fontsize{8}{10}\selectfont
\begin{tabu}{|c|[1pt]c@{\,}|c@{\,}|c@{\,}|c@{\,}|c@{\,}|c@{\,}|[1pt]c@{\,}|c@{\,}|c@{\,}|c@{\,}|c@{\,}|}\hline
i & 1 & 2 & 3 & \cdots & k & k+1 & k+2 & k+3 & \cdots & 2k & 2k+1 \\\tabucline[1pt]{-}
L=u_iw_i & 1 & 3 & 5 & \cdots & 2k-1 & 2k+1 & 2 & 4 & \cdots & 2k-2 & 2k \\\hline
R=v_iw_i & 3k+2 & 3k+1 & 3k & \cdots & 2k+3 & 2k+2 & 4k+2 & 4k+1 &  \cdots & 3k+4 & 3k+3 \\\hline
C_1=w_ix_{i,1}  &  6k+3 &  6k+2   & 6k+1 & \cdots & 5k+4 & 5k+3 & 5k+2 & 5k+1 & \cdots & 4k+4 & 4k+3      \\\hline
C_2=w_ix_{i,2}  &   10k+5 & 10k+4 & 10k+3 & \cdots & 9k+6 & 9k+5 & 9k+4 & 9k+3 & \cdots & 8k+6 & 8k+5     \\\hline
C_3=w_ix_{i,3}  &   6k+4 & 6k+5 & 6k+6 & \cdots & 7k+3 & 7k+4 & 7k+5 & 7k+6 & \cdots & 8k+3 & 8k+4     \\\hhline{|=*{10}{=}=|}
L_1=u_ix_{i,1}  &   14k+7 & 14k+5 & 14k+3 & \cdots & 12k+9 & 12k+7 & 14k+6 & 14k+4 & \cdots & 12k+10 & 12k+8       \\\hline
L_2=u_ix_{i,2}  &   18k+10 & 18k+12 & 18k+14 & \cdots & 20k+8 & 20k+10 & 18k+11 & 18k+13 & \cdots & 20k+7 & 20k+9        \\\hline
L_3=u_ix_{i,3}  &   18k+9 & 18k+7 & 18k+5 & \cdots & 16k+11 & 16k+9 & 18k+8 & 18k+6 & \cdots & 16k+12 & 16k+10       \\\hhline{|=*{10}{=}=|}
R_1=v_ix_{i,1}  &  22k+11 & 22k+9 & 22k+7 & \cdots & 20k+13 & 20k+11 & 22k+10 & 22k+8 & \cdots & 20k+14 & 20k+12     \\\hline
R_2=v_ix_{i,2}  &   11k+6 & 11k+7 & 11k+8 & \cdots & 12k+5 & 12k+6 & 10k+6 & 10k+7 & \cdots & 11k+4 & 11k+5      \\\hline
R_3=v_ix_{i,3}  &   14k+8 & 14k+10 & 14k+12 & \cdots & 16k+6 & 16k+8 & 14k+9 & 14k+11 & \cdots & 16k+5 & 16k+7       \\\hline \end{tabu}
\]

\nt Observe that
\begin{enumerate}[(a)]
  \item For each column, the sum of the first 5 row entries is a constant $25k+15$.
  \item For each column, the sum of the entries in rows $L, L_1, L_2, L_3$ (respectively, $R, R_1, R_2, R_3$) is a constant $50k+27$.
  \item For each $1\le a\le 3$, the sum of all the entries in rows $C_a, R_a, L_a$ is a constant $(2k+1)(39k+21)$.
\end{enumerate}

\begin{theorem}\label{thm-oddnP3VO3} For $k\ge 1$, $\chi_{la}((2k+1)P_3\vee O_3)=3$.\end{theorem}

\begin{proof} Clearly, $\chi_{la}((2k+1)P_3\vee O_3)\ge \chi((2k+1)P_3\vee O_3)=3$.

\nt Consider $(2k+1)$ copies of $P_3\vee O_3$ with the $i$-th copy vertex set $\{u_i, v_i, w_i, x_{i,1}, x_{i,2}, x_{i,3}\}$ and edge set $\{u_iw_i, v_iw_i, u_ix_{i,a}, v_ix_{i,a}, w_ix_{i,a}\;|\;1\le a\le 3\}$, $1\le i\le 2k+1$. Define a bijectively edge labeling $f : E((2k+1)(P_3\vee O_3)) \to [1, 11(2k+1)]$ such that for $1\le i\le 2k+1$ and $1\le a\le 3$:

\begin{enumerate}[(i)]
  \item $f(u_iw_i)$ is the $(L,i)$-entry,
  \item $f(v_iw_i)$ is the $(R,i)$-entry,
  \item $f(w_ix_{i,a})$ is the $(C_a,i)$-entry,
  \item $f(u_ix_{i,a})$ is the $(L_a,i)$-entry,
  \item $f(v_ix_{i,a})$ is the $(R_a,i)$-entry.
\end{enumerate}

\nt By the observations above, it is straightforward to conclude that $w_i$ has induced vertex label $25k+15$, $u_i$ and $v_i$ have same induced vertex label $50k+27$, and sum of the induced vertex label of $x_{i,1}$ (respectively, $x_{i,2}$ and $x_{i,3}$) for $1\le i\le 2k+1$ is $(2k+1)(39k+21)$. For $1\le a\le 3$, we now merge the vertices $x_{i,a}$, for $1\le i\le 2k+1$, as $x_a$. Then the resulting graph is $(2k+1)P_3\vee O_3$ and the induced vertex label for $x_a$ is $(2k+1)(39k+21)$ for each $a$. Thus, $f$ is a local antimagic 3-coloring and $\chi_{la}((2k+1)P_3\vee O_3)\le 3$. This completes the proof.
\end{proof}

\begin{example} For $k=1$, we have the following table

\[\fontsize{8}{10}\selectfont\begin{tabu}{|c|[1pt]c|c|c|}\hline
i & 1 & 2 & 3 \\\tabucline[1pt]{-}
L=u_iw_i & 1 & 3 & 2 \\\hline
R=v_iw_i & 5 & 4 & 6 \\\hline
C_1=w_ix_{i,1}  &  9 &  8 &7\\\hline
C_2=w_ix_{i,2}  &  15 & 14 & 13 \\\hline
C_3=w_ix_{i,3}  &  10 & 11 & 12 \\\hhline{|====|}
L_1=u_ix_{i,1}  &  21 & 19 & 20  \\\hline
L_2=u_ix_{i,2}  &  28 & 30 & 29 \\\hline
L_3=u_ix_{i,3}  &  27 & 25 & 26 \\\hhline{|====|}
R_1=v_ix_{i,1}  &  33 & 31 & 32  \\\hline
R_2=v_ix_{i,2}  &  17 & 18 & 16 \\\hline
R_3=v_ix_{i,3}  & 22 & 24 & 23 \\\hline
\end{tabu}\]
and get the following labeling for $3(P_3\vee O_3)$
\begin{figure}[H]
\centerline{\epsfig{file=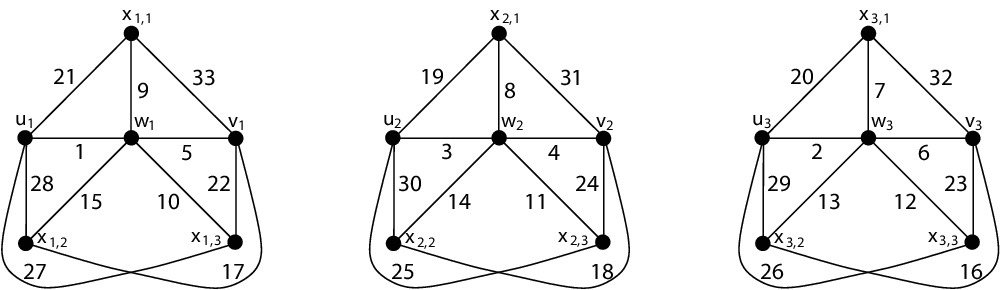, width=12cm}}
\caption{Graph $3(P_3\vee O_3)$.}\label{fig:3P3VO3}
\end{figure}
\nt For each $1\le a\le 3$, merge the vertices $x_{1,a}$, $x_{2,a}$ and $x_{3,a}$ to obtain the vertex $x_a$, we get $3P_3\vee O_3$ as required. \rsq
\end{example}

\section{Conclusion and Open problems}
In this paper, we constructed three $5\times (2k+1)$ matrices with entries in $[1,10k+5]$ that satisfy certain properties. Consequently, many tripartite graphs with local antimagic chromatic number 3 are obtained.

\ms\nt We shall in another paper construct two $5\times 2k$ matrices with entries in $[1,10k]$ to obtain many new bipartite or tripartite graphs with local antimagic chromatic number 3. This include showing that $\chi_{la}(FB(n)) = \chi_{la}(rFB(s)) = 3$ for even $n, s\ge 2$.

\ms\nt  The following problem arises naturally.

\begin{problem} Determine $\chi_{la}(G)$ for (i) $G= Pt(n)$, $n$ odd; (ii) $G= TB(n)$, $n$ odd; (iii) $G=rFB(s)$, $r$ even and $s$ odd.
\end{problem}

\nt Note that $FB(n) = nP_3 \vee K_1$ and that $\chi_{la}(nP_2\vee K_1)=3$ for $n\ge 1$, and $\chi_{la}(P_m \vee K_1) = 3$ for $m\ge 2$. (see~\cite{Arumugam, LSN, LScS}).

\begin{conjecture} For $m\ge 2$ and $n\ge 2$, $\chi_{la}(mP_n \vee K_1) = 3$. \end{conjecture}


\nt We end this paper with the following problem.

\begin{problem} Determine necessary and / or sufficient condition(s) for the existence of $m\times n$ matrices with entries in $[1,mn]$ that give rise to graph labelings and local antimagic chromatic numbers. \end{problem}

\end{document}